# Optimal oracle inequalities for model selection

**Charles Mitchell**

*e-mail:* `mitchell@stat.math.ethz.ch`
**and**

**Sara van de Geer**

*e-mail:* `geer@stat.math.ethz.ch`

**Abstract:** Model selection is often performed by empirical risk minimization. The quality of selection in a given situation can be assessed by risk bounds, which require assumptions both on the margin and the tails of the losses used. Starting with examples from the 3 basic estimation problems, regression, classification and density estimation, we formulate risk bounds for empirical risk minimization under successively weakening conditions and prove them at a very general level, for general margin and power tail behavior of the excess losses.

## 1. Introduction

Consider a sample $Z_1, \ldots, Z_N$ of independent random variables in some space $\mathcal{Z}$, whose distribution depends on an unknown parameter $f$. To estimate $f$, we split the sample into two parts: a test set $Z_1, \ldots, Z_n$ and a training set $Z_{n+1}, \ldots, Z_N$. Based on the training set various estimators of $f$ are constructed, say $\hat{f}_1, \ldots, \hat{f}_p$. To decide among these estimators, we use the test set. Suppose that $\gamma_f : \mathcal{Z} \to \mathbf{R}$ is a loss function. The final estimate $\hat{f}$ is now chosen to minimize the loss $\sum_{i=1}^{n} \gamma_{\hat{f}_j}(Z_i)$:

$$\hat{f} := \arg \min_{1 \le j \le p} \sum_{i=1}^{n} \gamma_{\hat{f}_j}(Z_i) \ .$$

In this note, we examine whether this procedure leads to taking, among the $p$ estimators, the "nearly best" one. Here, "nearly best" will be defined in terms of the excess risk of the estimators.

The behavior of the excess risk near the true value of $f$ will be called the margin behavior. We not only consider the classical case, which is quadratic margin behavior, but also more general margin behavior. For the tails of our excess loss functions, we consider both an exponential moment condition and a more general power tail condition. We prove a risk inequality under the most general combination of these conditions, and in doing so automatically obtain risk inequalities for more restricted situations. These latter situations represent examples we give from regression, classification and density estimation.

Note that the aggregation we perform is model selection aggregation. There is a rich body of literature on which aggregation method is best under a variety of conditions. Least-squares regression is discussed by (5), which gives the optimal rates of a number of methods, including linear and convex aggregation. A more general regression problem is addressed by (7). However, most of the literature deals with only one particular problem, such as regression, and also places strong conditions, like boundedness, on the functions and random variables involved. We obtain inequalities for a general loss function setup, and without boundedness conditions, at least when conditioning on the training set. Such conditioning on the training set is common practice; to average the results over the training data then requires more stringent conditions.

Another fairly general approach is found in (1), which looks at the general prediction problem, i.e. regression and classification, and uses a progressive mixture rule for aggregation, but with only a brief reference to averaging over the training stage, which would be part of the full sample splitting problem. On the other hand, (11) looks at sample splitting schemes with multiple splits and thus comes close to crossvalidation, but does so only for the problem of density estimation. A direct treatment of a crossvalidation scheme is to be found in (14). And in the context of classification, recent inequalities are given for recursive aggregation by mirror descent in (9) and for aggregation with exponential weights by (10).

### 1.1. Notation

The results will be conditional on the training set. We use $\mathbf{P}$ to denote the distribution of the test sample, and $\mathbf{E}$ denotes expectation of random variables depending on the test sample.

For $\gamma : \mathcal{Z} \to \mathbf{R}$, we write

$$P\gamma := \frac{1}{n} \sum_{i=1}^{n} \mathbf{E}\gamma(Z_i) \ ,$$

$$P_n\gamma := \frac{1}{n} \sum_{i=1}^{n} \gamma(Z_i) \ .$$

Let $\gamma_j : \mathcal{Z} \to \mathbf{R}$, $j = 1, \ldots, p$ be given loss functions in a class $\boldsymbol{\Gamma}$. Given the training set, $\gamma_j$ may be taken as short-hand (and slight abuse of) notation for $\gamma_{\hat{f}_j}$, $j = 1, \ldots, p$. We consider the estimator

$$\hat{\gamma} := \arg \min_{1 \le j \le p} P_n\gamma_j \ .$$

The target is

$$\gamma_0 := \arg \min_{\gamma \in \boldsymbol{\Gamma}} P\gamma \ .$$

The best approximation is

$$\gamma_* := \arg\min_{1 \le j \le p} P\gamma_j \ .$$

We define the excess risks

$$\hat{\mathcal{E}} := P(\hat{\gamma} - \gamma_0) \ ,$$

(which is a random variable),

$$\mathcal{E}_j := P(\gamma_j - \gamma_0)$$

and

$$\mathcal{E}_* := P(\gamma_* - \gamma_0) \ .$$

Without loss of generality, we assume that $\mathbf{\Gamma}$ is of the form $\mathbf{\Gamma} := \{\gamma_f : \ f \in \mathbf{F}\}$, where $\mathbf{F}$ is a subset of a metric space with metric $d$, and write (with some abuse of notation) $\gamma_{f_j}$ as $\gamma_j$, $\{f_j\}_{j=1}^p \subset \mathbf{F}$.

### *1.2. Goal*

Our goal is now to show that $\hat{\mathcal{E}}/\mathcal{E}_*$ is close to 1 (with large probability or in expectation). The results are modifications of inequalities of the form

$$(1-\delta)\mathbf{E}\hat{\mathcal{E}} \le (1+\delta)\mathcal{E}_* + \frac{\Delta_0}{\delta} \ ,$$

where $\delta > 0$ is an arbitrary small constant, and with $\Delta_0$ of order $\log(2p)/n$ and not depending on $\mathcal{E}_*$, see for example Chapter 7 in (6). In the standard setup of Section 4, we for instance show that for $1 \le m \le 1 + \log p$

$$\mathbf{E}\hat{\mathcal{E}}^{\frac{m}{2}} \le \left( \sqrt{\mathcal{E}_*} + \sqrt{\Delta_1} + \frac{\Delta_2}{\sqrt{\mathcal{E}_*}} \right)^m \ ,$$

with $\Delta_1$ and $\Delta_2$ both of order $\log(2p)/n$, and both not depending on $\mathcal{E}_*$. In particular, with $m = 2$, this reads

$$\mathbf{E}\hat{\mathcal{E}} \le \left( \sqrt{\mathcal{E}_*} + \sqrt{\Delta_1} + \frac{\Delta_2}{\sqrt{\mathcal{E}_*}} \right)^2 \ .$$

A sharp oracle inequality would be

$$\mathbf{E}\hat{\mathcal{E}} \le \mathcal{E}_* + \Delta_0 \ .$$

We conjecture that such sharpness cannot be established in a general setup by empirical risk minimization. Instead, e.g. mirror averaging could be used, see (8). See also (2) and (3) for some limitations of empirical risk minimization, and alternative approaches to overcome the limitations. We however believe empirical risk minimization remains an important topic of study because it is widely applied in practice, and is closely related to various cross validation schemes.

### *1.3. Convex loss*

In our proofs, we only use the property

$$P_n\hat{\gamma} \leq P_n\gamma_* \ .$$

In the convex case, this means sometimes that conditions can be weakened. Let **F** be a convex subset of a linear vector space, and suppose that $\mathbf{\Gamma} := \{\gamma_f : f \in \mathbf{F}\}$, with $f \mapsto \gamma_f$ convex, $P$-almost everywhere. Then for $0 \leq \alpha \leq 1$, we have the inequality

$$P_n\gamma_{\alpha\hat{f}+(1-\alpha)f_*} \leq \alpha P_n\hat{\gamma} + (1-\alpha)P_n\gamma_* \leq P_n\gamma_* \ .$$

This means that we can replace $\hat{\gamma}$ by $\gamma_{\alpha\hat{f}+(1-\alpha)f_*}$ throughout, leading to inequalities for the excess risk

$$\hat{\mathcal{E}}_\alpha = P\gamma_{\alpha\hat{f}+(1-\alpha)f_*} - P\gamma_0 \ .$$

From these, one may then often deduce inequalities for the original $d(\hat{f}, f_0)$. As we shall see, this extension (with $\alpha < 1$) allows us to work with weaker conditions (than with $\alpha = 1$). In particular, the example on maximum likelihood will use this approach with $\alpha$ set to $1/2$.

### *1.4. Organization of the paper*

The paper is organized as follows. Section 2 presents Bernstein's inequality. It is stated in the form of a probability inequality and a moment inequality. Section 3 presents the margin condition and some examples. In Section 4, we consider the standard setup with quadratic margin, and bounded loss or an exponential moment condition. Section 5 looks at loss with power moment conditions, and Section 6 at general margin behavior under the exponential moment condition, giving risk tail bounds. Section 7 formulates the general risk moment inequality, from which the previous specific results follow. Finally, the proofs are in Section 8.

## 2. Bernstein's inequality

Bernstein's inequality for a single average is well known, and the extension of Bernstein's probability inequality to a uniform probability inequality over $p$ averages is completely straightforward. The result can be seen as the simplest version of a concentration inequality in the spirit of e.g. (4) (emphasizing how tight these general concentration inequalities are). The moment inequality for the maximum of $p$ averages is perhaps less known.

For all $j$, we let

$$\gamma_j^c(\cdot) := \gamma_j(\cdot) - \mathbf{E}\gamma_j$$

denote the centered loss functions. To obtain our results, we we make assumptions on the tails of the centered excess losses $\gamma_j^c - \gamma_*^c$ or of their envelope $\Gamma := \max_{1 \le j \le p} |\gamma_j^c - \gamma_*^c|$ as follows:

**Definition 2.1.** *We say that the excess losses $\gamma_j - \gamma_*$ satisfy the exponential moment condition for some $K > 0$ if*

$$P\left|\gamma_j^c - \gamma_*^c\right|^m \le \frac{m!}{2}(2K)^{m-2}d^2(f_j, f_*) \tag{1}$$

*for all $m = 2, 3, \dots$ and for all $j = 1, \dots, p$.*

*We say that the envelope function $\Gamma$ has power tails of order $s > 1$ if there exists an $M \in (0, \infty)$ such that*

$$P(\{\Gamma > K\}) \le \left(\frac{M}{K}\right)^s \quad \forall K > 0 . \tag{2}$$

**Lemma 2.1.** *(Bernstein's inequality for the maximum of $p$ averages) Let loss functions $\gamma_j : \mathcal{Z} \to \mathbf{R}$, $j = 1, \dots, p$, be given. Assume that for some constant $K$ and for all $j$,*

$$P|\gamma_j^c|^m \le \frac{m!}{2}(2K)^{m-2}, \ m = 2, 3, \dots .$$

*Then for all $t > 0$,*

$$\mathbf{P}\left(\max_{1 \le j \le p} |P_n\gamma_j^c| \ge \sqrt{\frac{2(\log(2p) + t)}{n}} + \frac{2K(\log(2p) + t)}{n}\right) \le \exp[-t] . \tag{3}$$

*Moreover, for all $1 \le m \le 1 + \log p$,*

$$\left(\mathbf{E}\left(\max_{1 \le j \le p} |P_n\gamma_j^c|\right)^m\right)^{1/m} \le \sqrt{\frac{2\log(2p)}{n}} + \frac{2K\log(2p)}{n} . \tag{4}$$

In what follows, we will make repeated use of Bernstein's inequality. Hence, the term $2\log(2p)/n$ will appear frequently. From now on, we denote this term by

$$\Delta := \frac{2\log(2p)}{n} .$$

**Remark:** The moment inequality is for moments of order $m \leq 1 + \log p$. It can be extended to hold for general $m$, provided a slight adjustment, depending on $m$, is made on the constants. Because we have the situation in mind where $p$ is large, we have formulated the result for $m \leq 1 + \log p$ to facilitate the exposition.

**Corollary 2.1.** *(Weighted version of Bernstein's inequality) Assume that for some constant $K$, the condition*

$$P|\gamma_j^c - \gamma_*^c|^m \leq \frac{m!}{2}(2K)^{m-2}d^2(f_j, f_*), \; m = 2, 3, \ldots, \; \forall \; j \qquad (5)$$

*holds. Then for all $t > 0$ and $\tau > 0$*

$$\mathbf{P}\left(\max_{1 \leq j \leq p} \frac{|P_n(\gamma_j^c - \gamma_*^c)|}{d(f_j, f_*) \vee \tau} \geq \sqrt{\Delta + 2t/n} + \frac{K(\Delta + 2t/n)}{\tau}\right) \leq \exp[-t] \; .$$

*Moreover, for all $1 \leq m \leq 1 + \log p$,*

$$\left(\mathbf{E}\max_{1 \leq j \leq p}\left(\frac{|P_n(\gamma_j^c - \gamma_*^c)|}{d(f_j, f_*) \vee \tau}\right)^m\right)^{1/m} \leq \sqrt{\Delta} + \frac{K\Delta}{\tau} \; .$$

Define, for all $\gamma$, the variance

$$\sigma^2(\gamma) := P|\gamma^c|^2 \; .$$

Then clearly (5) implies that

$$d^2(f_j, f_*) \geq \sigma^2(\gamma_j - \gamma_*), \; \forall \; j \; .$$

Moreover, if the bound $|\gamma_j - \gamma_*| \leq 3K$ holds $\forall \; j$, then (5) holds with

$$d^2(f_j, f^*) = \sigma^2(\gamma_j - \gamma_*) \; \forall \; j.$$

In what follows, we will indeed often assume (5) with this value for $d(f_j, f^*)$, but we will also consider an extension. The choice of the metric $d$ is intertwined with the margin behavior, which we consider in the next section.

## 3. Margin behavior

**Definition 3.1.** *We say that the margin condition holds with strictly convex margin function $G(\cdot)$, if*

$$P(\gamma_j - \gamma_0) \geq G\left(d(f_j, f_0)\right), \ \forall \ j \ . \tag{6}$$

*Furthermore, we say that the margin condition holds with constants $\kappa > 1/2$ and $C > 0$, if (6) holds with*

$$G(u) = u^{2\kappa}/C^{2\kappa}, \ u > 0 \ .$$

As we shall see, $\kappa = 1$ in typical cases – but other, in particular larger, values can also occur.

Let us now consider some examples. In a regression or classification situation, we have i.i.d. random pairs $Z_i = (X_i, Y_i)$, with $Y_i \in \mathcal{Y} \subset \mathbf{R}$ a response variable, and $X_i \in \mathcal{X}$ a covariable, $i = 1, \ldots, n$. We then assume (for $i = 1, \ldots, n$) that the conditional distribution of $Y_i$, given $X_i = x$, only depends on $x$ and not on $i$. This can be done without loss of generality (as the index $i$ can be taken in as an additional covariable).

**Example 3.1.** *(Regression)* Suppose that $\{Z_i\}_{i=1}^n := \{(X_i, Y_i)\}_{i=1}^n$. Let $\mathbf{F}$ be a class of real-valued functions on $\mathcal{X}$, and for all $x \in \mathcal{X}$ and $y \in \mathcal{Y}$, let

$$\gamma_f(x, y) := \gamma(f(x), y), \ f \in \mathbf{F} \ .$$

Set

$$l(a, \cdot) = \mathbf{E}(\gamma(a, Y_i)|X_i = \cdot), \ a \in \mathbf{R} \ .$$

We moreover write $l_f(x) := l(f(x), x)$. As target we take the overall minimizer

$$f_0(\cdot) := \arg\min_{a \in \mathbf{R}} l(a, \cdot) \ .$$

We now check whether the margin condition holds with $\kappa = 1$ and

$$d^2(f, f_0) := K_2^2 P|f - f_0|^2 \ ,$$

where $K_2$ is an appropriate constant.

**Lemma 3.1.** *Assume that for some $K_1 > 0$, and all $|f - f_0| \leq K_1$,*

$$l_f - l_{f_0} \geq (f - f_0)^2/C_0^2 \ , \tag{7}$$

*Then*

$$P(\gamma_f - \gamma_{f_0}) \geq d^2(f, f_0)/C^2 \ ,$$

*where $C^2 := C_0^2 K_2^2$. If we moreover assume (for $i = 1, \ldots, n$) that*

$$\text{var}(\gamma_f(Z_i) - \gamma_{f_0}(Z_i)) \leq K_2^2 \mathbf{E}(f(X_i) - f_0(X_i))^2 \ , \tag{8}$$

*then for all $\|f - f_0\|_\infty \leq K_1$, we have*

$$\sigma^2(\gamma_f - \gamma_{f_0}) \leq d^2(f, f_0) \ .$$

If $l(a, \cdot)$ has two derivatives near $a = f_0(\cdot)$, and the second derivatives are positive and bounded away from zero, then $l(a, \cdot)$ behaves quadratically near its minimum, i.e., then (7) holds for some $K_1 > 0$.

It also also clear that (8) holds as soon as $\gamma(\cdot, y)$ is Lipschitz for all $y$, with Lipschitz constant $L$. Then we may take $K_2 = L$. When $\gamma(\cdot, y)$ is not Lipschitz (e.g., quadratic loss), it may be useful to define

$$e_f(Z_i) := \gamma(f(X_i), Y_i) - l_f(X_i) .$$

Then obviously

$$\text{var}(\gamma_f(Z_i) - \gamma_{f_0}(Z_i)) = \text{var}(e_f(Z_i) - e_{f_0}(Z_i)) + \text{var}(l_f(X_i) - l_{f_0}(X_i)) . \quad (9)$$

Note that with fixed design, the second term in (9) vanishes.

*Quadratic loss:*

In the case of least squares, the loss function is

$$\gamma(f, y) := (y - f)^2 ,$$

Then

$$l_f - l_{f_0} = |f - f_0|^2 ,$$

and

$$e_f(Z_i) - e_{f_0}(Z_i) = 2\epsilon_i(f(X_i) - f_0(X_i)) ,$$

with $\epsilon_i := Y_i - f_0(X_i)$. Assuming that the conditional variance is bounded by some constant $\sigma_\epsilon$, i.e.,

$$\max_{1 \le i \le n} \text{var}(Y_i | X_i) \le \sigma_\epsilon^2 , \quad (10)$$

we may conclude the following.

*Least squares with fixed design:*

The margin condition holds with $\kappa = 1$ and $C^2 = 4\sigma_\epsilon^2$ .

*Least squares with random design:*

If $\|f_j - f_0\|_\infty \le K_1$ for all $j$, the margin condition holds with $\kappa = 1$ and $C^2 = 4\sigma_\epsilon^2 + K_1^2$ .

**Example 3.2.** *(Classification)* Suppose that $Z_i = (X_i, Y_i)$, with $Y_i \in \mathcal{Y} := \{0, 1\}$ a label, $i = 1, \dots, n$. Let $\mathbf{F}$ be a class of functions $f : \mathcal{X} \to [0, 1]$. We consider 0/1-loss

$$\gamma_f(x, y) = \gamma(f(x), y) := (1 - y)f(x) + y(1 - f(x)), \ f \in \mathbf{F}, \ (x, y) \in \mathcal{X} \times \{0, 1\} .$$

For $a \in [0, 1]$, write

$$l(a, \cdot) := \mathbf{E}(\gamma(a, Y_i) | X_i = \cdot)$$

$$= (1 - \eta)a + \eta(1 - a) = a(1 - 2\eta) + \eta ,$$

where $\eta = \mathbf{E}(Y_i | X_i = \cdot)$. The target is again the overall minimizer

$$f_0 := \arg \min_{a \in [0, 1]} l(a, \cdot) .$$

It is clear that $f_0$ is the Bayes rule

$$f_0 = \mathrm{l}\{1 - 2\eta < 0\} + q\{1 - 2\eta = 0\} \ ,$$

with $q$ an arbitrary value in $[0, 1]$. We moreover have

$$P(\gamma_f - \gamma_{f_0}) = P|(f - f_0)(1 - 2\eta)| \ .$$

Consider the functions

$$H_1(v) \leq vP\mathrm{l}\{|1 - 2\eta| < v\}, \ v \in [0, 1] \ ,$$

and

$$G_1(u) = \max_v\{uv - H_1(v)\}, \ u \in [0, 1]$$

(assuming the maximum exists).

**Lemma 3.2.** *The inequality*

$$P(\gamma_f - \gamma_{f_0}) \geq G\bigg(\sigma(\gamma_f - \gamma_{f_0})\bigg)$$

*holds with $G(u) = G_1(u^2)$, $u \in [0, 1]$.*

If $H_1(v) = 0$ for $v \leq C_1$, we take $G_1(u) = C_1 u$. More generally, the Tsybakov margin condition (see (12)) assumes that one may take, for some $C_1 \geq 1$ and $\gamma \geq 0$,

$$H_1(v) = v(C_1 v)^{1/\gamma} \ ,$$

Then one has

$$G_1(u) = u^{1+\gamma}/C^{1+\gamma}$$

where

$$C = C_1^{\frac{1}{1+\gamma}} \gamma^{-\frac{\gamma}{1+\gamma}} (1 + \gamma) \ .$$

Thus, then the margin condition holds with this value of $C$ and with $\kappa = 1 + \gamma$ (and for any $d$ satisfying $d(f_j, f_0) \geq \sigma(\gamma_j - \gamma_0)$, $\forall j$).

**Example 3.3.** *(Maximum likelihood)* Suppose that $\{Z_i\}_{i=1}^n$ are iid. with density $f_0 := dP/d\mu$, where $\mu$ is a $\sigma$-finite dominating measure. Let $\mathbf{F}$ be a (convex, say) class of densities w.r.t. $\mu$, containing $f_0$. Consider the transformed log-likelihood loss

$$\gamma_f(\cdot) := \gamma(f(\cdot)),$$

where $\gamma(a) = -\log(a)/2$. Define

$$\bar{f} = (f + f_*)/2, \ f \in \mathbf{F} \ .$$

The squared Hellinger distance of densities $f$ and $\tilde{f}$ is

$$h^2(f, \tilde{f}) = \frac{1}{2} \int \left(\sqrt{f} - \sqrt{\tilde{f}}\right)^2 d\mu, \ f, \tilde{f} \in \mathbf{F} \ .$$

We now check the margin condition with $\kappa = 1$ and $d(f, f_0) = Ch(f, f_0)$.

**Lemma 3.3.** *For all densities $f$, we have*

$$P(\gamma_f - \gamma_{f_0}) \geq h^2(f, f_0) \ .$$

*Moreover, under the assumption*

$$\sqrt{\frac{f_0}{f_*}} \leq \frac{C}{8} \ ,$$

*we have*

$$\sigma(\gamma_{\bar{f}} - \gamma_{f_*}) \leq Ch(\bar{f}, f_*) \ .$$

## 4. Quadratic margin and exponential moments

The first case we shall look at is the one with quadratic margin condition ($\kappa = 1$) and exponential moments on the tails of the loss functions. This encompasses e.g. regression with sub-Gaussian errors, as well as many situations where estimators and losses have absolute bounds.

### *4.1. General loss*

**Lemma 4.1.** *Suppose that the margin condition holds, with constants $\kappa = 1$ and $C > 0$. Assume moreover that the loss functions satisfy the exponential moment condition (1) for some $K > 0$. Then for all $t > 0$, and when $\mathcal{E}_* > 0$,*

$$\mathbf{P}\left(\sqrt{\hat{\mathcal{E}}} \geq \sqrt{\mathcal{E}_*} + C\sqrt{\Delta + 2t/n} + \frac{K(\Delta + 2t/n)}{\sqrt{\mathcal{E}_*}}\right) \leq e^{-t} \ .$$

*When*

$$\mathcal{E}_* \leq K(\Delta + 2t/n) \ ,$$

*we have*

$$\mathbf{P}\left(\sqrt{\hat{\mathcal{E}}} \geq (C + 2\sqrt{K})\sqrt{\Delta + 2t/n}\right) \leq e^{-t} \ .$$

*Moreover, for all $1 \leq m \leq 1 + \log p$, when $\mathcal{E}_* > 0$, we have*

$$\left\|\sqrt{\frac{\hat{\mathcal{E}}}{\mathcal{E}_*}}\right\|_m \leq 1 + C\sqrt{\frac{\Delta}{\mathcal{E}_*}} + \frac{K\Delta}{\mathcal{E}_*} \ ,$$

*and when*

$$\mathcal{E}_* \leq K\Delta$$

*we have*

$$\left\|\sqrt{\hat{\mathcal{E}}}\right\|_m \leq (C + 2\sqrt{K})\sqrt{\Delta} \ .$$

*Proof.* All statements in this lemma are special cases of Lemma 6.2. □

**Corollary 4.1.** *(Asymptotics) When*

$$\mathcal{E}_* \gg \left(K + C^2\right)\Delta \ ,$$

*it holds (for $m \leq 1 + \log p$) that*

$$\mathbf{E}\left|\sqrt{\frac{\hat{\mathcal{E}}}{\mathcal{E}_*}}\right|^m \to 1 \ .$$

### *4.2. Maximum likelihood*

Define

$$\hat{\mathcal{K}} := P(\gamma_{(\hat{f}+f_*)/2} - \gamma_{f_0}) = \hat{\mathcal{E}}_{1/2} \ ,$$

and

$$\mathcal{K}_* := P(\gamma_{f_*} - \gamma_{f_0}) = \mathcal{E}_* \ .$$

Note that $\hat{\mathcal{K}}$ and $\mathcal{K}_*$ are Kullback-Leibler information numbers. Lemma 4.2 below presents a version of Lemma 4.1 for the maximum likelihood framework.

**Lemma 4.2.** *Suppose that*

$$\sqrt{\frac{f_0}{f_*}} \le \frac{C}{8} \ .$$

*Then for all $t > 0$, and when $\mathcal{K}_* > 0$,*

$$\mathbf{P}\left( \sqrt{\frac{\hat{\mathcal{K}}}{\mathcal{K}_*}} \ge 1 + C\sqrt{\frac{\Delta + 2t/n}{\mathcal{K}_*}} + \frac{\Delta + 2t/n}{\mathcal{K}_*} \right) \le \mathrm{e}^{-t} \ .$$

*When*

$$\mathcal{K}_* \le \Delta + 2t/n \ ,$$

*we have*

$$\mathbf{P}\left( \sqrt{\hat{\mathcal{K}}} \ge (C+2)\sqrt{\Delta + 2t/n} \right) \le \mathrm{e}^{-t} \ ,$$

*Moreover, for all $1 \le m \le 1 + \log p$, when $\mathcal{K}_* > 0$, we have*

$$\left\| \sqrt{\frac{\hat{\mathcal{K}}}{\mathcal{K}_*}} \right\|_m \le 1 + C\sqrt{\frac{\Delta}{\mathcal{K}_*}} + \frac{\Delta}{\mathcal{K}_*} \ ,$$

*and when*

$$\mathcal{K}_* \le \Delta \ ,$$

*we have*

$$\left\| \sqrt{\hat{\mathcal{K}}} \right\|_m \le (C+2)\sqrt{\Delta} \ .$$

## 5. Quadratic margin, power tails

### *5.1. Large values of p, power tails of the envelope function*

**Lemma 5.1.** *Suppose that the margin condition holds, with constants $\kappa = 1$ and $C > 0$, and some $d$ satisfying $d(f_j, f_0) \ge \sigma(\gamma_j - \gamma_0)$, $\forall\ j$. Assume moreover that the envelope has power tails, i.e., that (2) holds for some $s > 1$ and $M \in (0, \infty)$. Then for $1 \le m \le 1 + \log p$, and $m < 2s$, when $\mathcal{E}_* > 0$,*

$$\left\| \sqrt{\frac{\hat{\mathcal{E}}}{\mathcal{E}_*}} \right\|_m \le 1 + C\sqrt{\frac{\Delta}{\mathcal{E}_*}} + c_{m,s}\left( \frac{M}{\mathcal{E}_*} \right)^{\frac{2s}{2s+m}} \Delta^{\frac{2s-m}{2s+m}} \ ,$$

*where*

$$c_{m,s} := \left(\frac{m}{2s-m}\right)^{\frac{2}{2s+m}} \mathcal{C}\left(\frac{2s-m}{2m}\right) \ ,$$

*with*

$$\mathcal{C}(a) = a^{\frac{1}{1+a}} + a^{-\frac{a}{1+a}}, \ a > 0 \ .$$

*Moreover, if*

$$\mathcal{E}_* \le c_{m,s}^{\frac{2s+m}{2s}} \left(\frac{2s-m}{2s+m}\right)^{\frac{2s+m}{2s}} M\Delta^{\frac{2s-m}{2s}} \ ,$$

*we have*

$$\left\|\sqrt{\hat{\mathcal{E}}}\right\|_m \le C\sqrt{\Delta} + c'_{m,s}\sqrt{M}\Delta^{\frac{2s-m}{4s}} \ ,$$

*where*

$$c'_{m,s} := c_{m,s}^{\frac{2s+m}{4s}} \mathcal{C}\left(\frac{2s-m}{2s+m}\right) \ .$$

*Proof.* The first moment inequality is a special case of Theorem 7.1(i), and the other statements are immediate consequences of it.    □

**Corollary 5.1.** *(Asymptotics) When*

$$\mathcal{E}_* \gg C^2\Delta + c_{2,s}^{\frac{s+1}{s}} M\Delta^{\frac{s-1}{s}} \ ,$$

*then we have*

$$\mathbf{E}\left(\frac{\hat{\mathcal{E}}}{\mathcal{E}_*}\right) \to 1 \ .$$

### *5.2. Lower bounds*

#### *5.2.1. Large values of $p$*

Section 5.2.2 will show that (with $m = 2$) Lemma 5.1 can be improved if $p$ is small compared to $\sqrt{n}$. In this section, we present a lower bound where $p = \sqrt{n} + 1$ (or larger), which shows that essentially, Lemma 5.1 (with $m = 2$) cannot be improved. For a fair comparison, the same conditions are imposed as in Lemma 4.1: the margin condition, and the tail condition.

We consider quadratic loss

$$\gamma_f(\cdot, y) = (y - f)^2 \ .$$

Morover, we let $X_1, \ldots, X_n$ be fixed and

$$Y_i = f_0(X_i) + \epsilon_i, \ i = 1, \ldots, n \ ,$$

where $\epsilon_1, \ldots, \epsilon_n$ are i.i.d. copies of a random variable $\epsilon$, which has a double Pareto distribution, with parameter $s > 2$, i.e., the distribution of $\epsilon$ is symmetric around 0, and

$$P(|\epsilon| \le u) = 1 - \frac{1}{(1+u)^s}, \ u > 0 \ .$$

Now, suppose $p = \sqrt{n} + 1$, $f_p := f_0 \equiv 0$, and that for $j = 1, \ldots, p - 1 = \sqrt{n}$,

$$f_j(x) = 1\{x = X_j\} n^{\frac{1}{2s}}, \ x \in \mathcal{X} \ .$$

**Lemma 5.2.** *The margin condition holds with $\kappa = 1$ and $C^2 = 8/((s-2)(s-1))$, and the power tail condition (2) holds with $M = 2$. Moreover, for $n \geq 2^{2s}$, with probability at least $1 - \exp[-2^{-s}]$ we have*

$$\hat{\mathcal{E}} \geq n^{-\frac{s-1}{s}} \ .$$

**Remark** One may easily extend the situation to $p \gg \sqrt{n}$, because one may add, as candidates, as many bounded functions $f_j$, say $\|f_j\|_\infty \leq 1$, without destroying the moment condition (increasing $M$ from $M = 2$ to $M = 4$). These added functions may be selected by the least squares estimator, but if they all all have norm $P f_j^2 \geq n^{-\frac{s-1}{s}}$, selecting one of those still gives the same lower bound.

*5.2.2. Small values of p: the least squares case*

We consider again quadratic loss, and

$$Y_i = f_0(X_i) + \epsilon_i, \ i = 1, \ldots, n \ ,$$

with fixed design $X_1, \ldots, X_n$ and $\epsilon_1, \ldots, \epsilon_n$ are i.i.d. copies of a random variable $\epsilon$ with mean zero. Assume now a finite $s$-th moment

$$M^s := E|\epsilon|^s \ .$$

We now show that a lower bound of order $n^{-\frac{s-1}{s}}$ for $\mathbf{E}\hat{\mathcal{E}}$ will not hold if $p$ is small compared to $\sqrt{n}$.

**Lemma 5.3.** *We have*

$$\left( \mathbf{E} \left( \sqrt{\hat{\mathcal{E}}} \right)^s \right)^{\frac{1}{s}} \leq C c_s p^{1/s} M / \sqrt{n} + \sqrt{\mathcal{E}_*} \ ,$$

*where*

$$c_s := 2\sqrt{\frac{2}{\pi}} \Gamma^{1/s} \left( \frac{s+1}{2} \right) \ .$$

**Corollary 5.2.** *If $p \leq \sqrt{n}$ it holds that*

$$\mathbf{E}\hat{\mathcal{E}} \leq \left( C c_s n^{-\frac{s-1}{2s}} M + \sqrt{\mathcal{E}_*} \right)^2 \ .$$

## 6. General margin, exponential moments

In this section, we weaken the margin condition to allow for parameter values $\kappa > 1$. Example 3.2 already showed us the necessity of this more general condition, as it overlaps with Tsybakov's margin condition.

**Lemma 6.1.** *Suppose that the margin condition holds, with strictly convex margin function $G$. Let $H$ be the convex conjugate of $G$. Assume that for some $r \leq 1 + \log p$, the function $H(v^{\frac{1}{r}})$, $v > 0$, is concave. Assume moreover that the exponential moment condition (1) holds for some $K > 0$. Then for all $0 < \delta < 1$, and $\varepsilon > 0$, we have*

$$(1 - \delta)\mathbf{E}\hat{\mathcal{E}} \leq 2\delta H\left(\frac{\sqrt{\Delta}}{\delta} + \frac{K\Delta}{2\delta G^{-1}(\mathcal{E}_* \vee \varepsilon)}\right) + (1 + \delta)\mathcal{E}_* .$$

Lemma 7.1 is already set in the form of a non-sharp oracle inequality, rather than as a general bound on risk moments. In Section 7, we will derive a similar oracle inequality for the margin condition with $G(u) = u^{2\kappa}/C^{2\kappa}$, but first we give the more general risk bound in this case:

**Lemma 6.2.** *Suppose that the margin condition holds, with constants $\kappa \geq 1$ and $C > 0$. Assume moreover that the exponential moment condition (1) holds for $K > 0$. Then for $\mathcal{E}_* > 0$, and all $t > 0$, we have*

$$\mathbf{P}\left(\hat{\mathcal{E}}^{\frac{1}{2\kappa}} \geq \mathcal{E}_*^{\frac{1}{2\kappa}} + A(\kappa)\left(C\sqrt{\Delta + 2t/n} + \frac{K(\Delta + 2t/n)}{\mathcal{E}_*^{\frac{1}{2\kappa}}}\right)^{\frac{1}{2\kappa - 1}}\right) \leq \mathrm{e}^{-t} ,$$

*where*

$$A(\kappa) := \frac{1 + (2\kappa - 1)^{\frac{1}{2\kappa - 1}}}{(2\kappa)^{\frac{1}{2\kappa - 1}}} < 2 .$$

*Moreover, for*

$$\mathcal{E}_* \leq K(\Delta + 2t/n) ,$$

*we have*

$$\mathbf{P}\left(\hat{\mathcal{E}}^{\frac{1}{2\kappa}} \geq A(\kappa)\left(C\sqrt{\Delta + 2t/n}\right)^{\frac{1}{2\kappa - 1}} + 2\left(K(\Delta + 2t/n)\right)^{\frac{1}{2\kappa}}\right) \leq \mathrm{e}^{-t} .$$

*Furthermore, for all $m \leq (1 + \log p)(2\kappa - 1)$, when $\mathcal{E}_* > 0$,*

$$\left\|\left(\hat{\mathcal{E}}\right)^{\frac{1}{2\kappa}}\right\|_m \leq \mathcal{E}_*^{\frac{1}{2\kappa}} + A(\kappa)\left(C\sqrt{\Delta} + \frac{K\Delta}{\mathcal{E}_*^{\frac{1}{2\kappa}}}\right)^{\frac{1}{2\kappa - 1}} ,$$

*and when*

$$\mathcal{E}_* \leq K\Delta ,$$

*we have*

$$\left\|\hat{\mathcal{E}}^{\frac{1}{2\kappa}}\right\|_m \leq A(\kappa)\left(C\sqrt{\Delta}\right)^{\frac{1}{2\kappa - 1}} + 2\left(K\Delta\right)^{\frac{1}{2\kappa}} .$$

*Proof.* The moment inequalities follow from Theorem 7.1(ii), first taking $\tau^2 := \mathcal{E}_*$, and then $\tau^2 := K\Delta$. The tail bounds follow from the same theorem by taking first $\tau^2 := \mathcal{E}_*$, then $\tau^2 := K(\Delta + 2t/n)$. $\qquad\square$

**Corollary 6.1** (Asymptotics)**.** *When*

$$\mathcal{E}_* \gg C\Delta^{\frac{\kappa}{2\kappa-1}} + K\Delta \ ,$$

*it holds (for all $m \leq (1 + \log p)(2\kappa - 1)$) that*

$$\mathbf{E}\left(\frac{\hat{\mathcal{E}}}{\mathcal{E}_*}\right)^{\frac{m}{\kappa}} \to 1 \ .$$

## 7. General margin & tails

We now formulate our main theorem, whose proof also contains the proof of the moment bounds in Lemma 6.2:

**Theorem 7.1.**     *(i) Suppose that the margin condition holds for the loss functions $\gamma_j$ with constants $\kappa \geq 1$ and $C > 0$ and some $d$ satisfying $d(f_j, f_0) \geq \sigma(\gamma_j - \gamma_0)$, $\forall\ j$. Also assume that the envelope $\Gamma$ has power tails in the form of (2) for some $s > 1$ and $M > 0$. Then for all $m$ in the interval $[2\kappa, \min(2s\kappa, 1 + \log(p))[$ and for all $\tau > 0$, we have the following inequality:*

$$\left\|\left(\hat{\mathcal{E}}\right)^{\frac{1}{2\kappa}}\right\|_m \leq (\mathcal{E}_* \vee \tau)^{\frac{1}{2\kappa}} + A(\kappa) \cdot C^\alpha \cdot \Delta^{\alpha/2}$$

$$+\xi(\kappa, s, m) \cdot M^{\frac{s}{m}} \cdot \frac{\alpha}{\alpha+\beta} \cdot \Delta^{\frac{\alpha\beta}{\alpha+\beta}} \cdot (\mathcal{E}_* \vee \tau)^{-\frac{1}{2\kappa} \cdot \frac{\alpha\beta}{\alpha+\beta}} \ ,$$

*where*

$$\alpha := \frac{1}{2\kappa-1} \ , \quad \beta := \frac{s}{m} - \frac{1}{2\kappa} \ ,$$

$$A(\kappa) := \frac{1 + (2\kappa-1)^{\frac{1}{2\kappa-1}}}{\kappa^{\frac{1}{2\kappa-1}}}$$

*and*

$$\xi(\kappa, s, m) := A(\kappa)^{\frac{\beta}{\alpha+\beta}} \cdot 2^{\frac{1}{2\kappa} \cdot \frac{\alpha}{\alpha+\beta}} \cdot \left(\frac{m}{2s\kappa-m}\right)^{\frac{\alpha}{\alpha+\beta} \cdot \frac{1}{m}} \cdot \left(\left(\frac{\beta}{\alpha}\right)^{\frac{\alpha}{\alpha+\beta}} + \left(\frac{\alpha}{\beta}\right)^{\frac{\beta}{\alpha+\beta}}\right) \ .$$

*(ii) Furthermore, if the excess losses satisfy the exponential moment condition (1) for some constants $K > 0$, then*

$$\left\|\left(\hat{\mathcal{E}}\right)^{\frac{1}{2\kappa}}\right\|_m \leq (\mathcal{E}_* \vee \tau)^{\frac{1}{2\kappa}} + A(\kappa) \cdot \left(C \cdot \sqrt{\Delta} + \frac{K\Delta}{(\mathcal{E}_* \vee \tau)^{\frac{1}{2\kappa}}}\right)^{\frac{1}{2\kappa-1}} \ .$$

*In this case we also have tail bounds*

$$\mathbf{P}\left(\hat{\mathcal{E}}^{\frac{1}{2\kappa}} \geq \mathcal{E}_*^{\frac{1}{2\kappa}} + A(\kappa)\left(C\sqrt{\Delta + 2t/n} + \frac{K(\Delta + 2t/n)}{\mathcal{E}_*^{\frac{1}{2\kappa}}}\right)^{\frac{1}{2\kappa-1}}\right) \leq \mathrm{e}^{-t}$$

*for all $t > 0$ .*

These statements lead to simpler ones if we use that $\tau \leq \mathcal{E} \vee \tau \leq \mathcal{E} + \tau$ and then optimize over $\tau$:

**Corollary 7.1.** *Under the conditions of Theorem 7.1, we have the inequality*

$$(i) \qquad \left\|\left(\hat{\mathcal{E}}\right)^{\frac{1}{2\kappa}}\right\|_m \leq \mathcal{E}_*^{\frac{1}{2\kappa}} + A(\kappa)\cdot C^\alpha \cdot \Delta^{\alpha/2} + \tilde{\xi}(\kappa,s,m)\cdot M^{\frac{s}{m}\cdot\frac{\alpha}{\alpha+\beta+\alpha\beta}}\cdot\Delta^{\frac{\alpha\beta}{\alpha+\beta+\alpha\beta}}$$

*when the loss envelope $\Gamma$ has power tails (2), and*

$$(ii) \qquad \left\|\left(\hat{\mathcal{E}}\right)^{\frac{1}{2\kappa}}\right\|_m \leq \mathcal{E}_*^{\frac{1}{2\kappa}} + \cdot C^{\frac{1}{2\kappa-1}}\cdot\Delta^{\frac{1}{4\kappa-2}} + \sqrt{A(\kappa)}\cdot K^{\frac{1}{2\kappa}}\cdot\Delta^{\frac{1}{2\kappa}} \ .$$

*when the excess losses satisfy the exponential moment condition (1).*

### 7.1. Special cases of Corollary 7.1

We can apply Corollary 7.1 to the (more restricted) cases described in the previous sections:

**Quadratic margin, power tails:**  Here $\kappa = 1$ and thus $\alpha = 1$, $\beta = s/m - 1/2$ and $A(\kappa) = 1$. Theorem 7.1 thus implies

$$\left\|\sqrt{\hat{\mathcal{E}}}\right\|_m \leq \sqrt{\mathcal{E}_*} + \cdot C\cdot\sqrt{\Delta} + \xi(1,s,m)\cdot M^{\frac{2s}{2s+m}}\cdot\Delta^{\frac{2s-m}{2s+m}}\cdot(\mathcal{E}_*)^{\frac{m-2s}{2(2s+m)}} \ ,$$

as in Lemma 5.1; the corresponding simplified version from Corollary 7.1 is

$$\left\|\sqrt{\hat{\mathcal{E}}}\right\|_m \leq \sqrt{\mathcal{E}_*} + C\cdot\sqrt{\Delta} + \xi(1,s,m)\cdot\sqrt{M}\cdot\Delta^{1-\frac{1}{s}} \ .$$

For $m = 2$, this implies

$$P\hat{\mathcal{E}} \leq (1+\delta)\sqrt{\mathcal{E}_*} + \frac{1}{\delta}\cdot\left(C\cdot\Delta + \xi(1,s,2)\cdot M\cdot\Delta^{\frac{s-1}{2s}}\right) \ .$$

In the example of least-squares regression (Example 3.1), we know that a quadratic margin condition holds, e.g. for the fixed design with $C^2 := 4\sigma_\epsilon^2$. If furthermore we assume that the errors $\epsilon_i$ possess some finite moment of order $2s > 2$ – a less restrictive assumption than the Gaussianity often assumed – then the loss has power tails of order $s > 1$:

$$\gamma_f(x,y) = \gamma(f(x),y) = (y-f(x))^2 = (\epsilon + f_0(x) - f(x))^2$$

$$\Rightarrow \mathbf{E}\left[\Gamma^s\right] \leq 2^s \mathbf{E}\left[\sup_{f\in F}\left|\gamma_f^c(X,Y)\right|^s\right] \leq 2^{4s-1} \cdot \mathbf{E}\left[|\epsilon|^{2s} + \sup_{f\in F}|f_0(X)-f(X)|^{2s}\right]$$

$$= 2^{4s-1} \cdot \left(\mathbf{E}|\epsilon|^{2s} + \mathbf{E}\sup_{f\in F}|f_0(X)-f(X)|^{2s}\right) =: M \ ,$$

and so by Chebyshev,

$$P\left(\{\Gamma > K\}\right) \leq \frac{\mathbf{E}\left[\Gamma^s\right]}{K^s} \leq \left(\frac{M}{K}\right)^s \quad \forall K > 0 \ .$$

**General margin, exponential tails**   The risk bound in this case was given in Part (ii) of Corollary 7.1, whose correction term is of order $O(\Delta^{1/(4\kappa-2)})$. This leads to an oracle inequality of the form

$$P\hat{\mathcal{E}} \leq (1+\delta)\mathcal{E}_* + \frac{1}{\delta}O\left(\Delta^{\frac{2\kappa}{4\kappa-2}}\right) \quad \forall \delta > 0 \ .$$

In Example 3.2, we have already seen the margin condition for $C = C_1^{1/(1+\gamma)}\gamma^{-\gamma/(1+\gamma)}(1+\gamma)$ and $\kappa = 1+\gamma$, where $\gamma \geq 0$, as a consequence of Tsybakov's margin condition. Furthermore,

$$\begin{aligned}
P\left|\gamma_f^c - \gamma_{f_0}^c\right|^m &= P\left|(f(X)-f_0(X))\cdot(1-2Y) - P\left|(f-f_0)(1-2\eta)\right|\right|^m \\
&\leq 2^{m-2}\cdot P\left|\gamma_f^c - \gamma_{f_0}^c\right|^2 = 2^{m-2}\cdot \sigma^2\left(\gamma_f - \gamma_{f_0}\right)
\end{aligned}$$

for all $f$ in this example, which means that the excess losses have exponential moments (1) with $K = 1$. Thus we have an oracle inequality

$$\begin{aligned}
\hat{\mathcal{E}} &\leq (1+\delta)\mathcal{E}_* + \frac{1}{\delta}\left(\tilde{A}_1(C_1,\gamma)\cdot\Delta^{\frac{1+\gamma}{1+2\gamma}} + \tilde{A}_2\cdot\Delta\right) \\
&= \mathcal{E}_*^{\frac{1}{2\kappa}} + O(\Delta^{\frac{1+\gamma}{1+2\gamma}}) \ .
\end{aligned}$$

## 8. Proofs

### *8.1. Proofs for Section 2*

**Proof of Lemma 2.1.** Without loss of generality, suppose that $\mathbf{E}\gamma_j(Z_i) = 0$ for all $i$ and $j$. Bernstein's probability inequality says that for all $t > 0$,

$$\mathbf{P}\left(\frac{1}{n}\sum_{i=1}^{n}\gamma_j(Z_i) \geq 2Kt + \sqrt{2t}\right) \leq \exp\left[-nt\right], \forall\ j\ . \tag{11}$$

This inequality follows from the intermediate result

$$\mathbf{E}\exp\left[\sum_{i=1}^{n}\gamma_j(Z_i)/L\right] \leq \exp\left[\frac{n}{2(L^2 - 2LK)}\right],\ \forall\ j\ , \tag{12}$$

which holds for all $L > 2K$. Inequality (3) follows immediately from (11).

To prove (4), we apply Lemma 8.1. We then obtain for all $L > 0$, and all $m$

$$\mathbf{E}\left(\max_j |\sum_{i=1}^{n}\gamma_j(Z_i)|^m\right) \leq L^m \log^m\left[\mathbf{E}\exp[\max_j |\sum_{i=1}^{n}\gamma_j(Z_i)|/L] - 1 + \mathrm{e}^{m-1}\right].$$

From (12), and invoking $\mathrm{e}^{|x|} \leq \mathrm{e}^{x} + \mathrm{e}^{-x}$, we obtain for $L > 2K$,

$$L^m \log^m\left[\mathbf{E}\exp[\max_j |\sum_{i=1}^{n}\gamma_j(Z_i)|/L] - 1 + \mathrm{e}^{m-1}\right]$$

$$\leq L^m \log^m\left[p\{2\exp\left[\frac{n}{2(L^2 - 2LK)}\right] - 1\} + \mathrm{e}^{m-1}\right]$$

$$\leq L^m \log^m\left[(2p + \mathrm{e}^{m-1} - p)\exp\left[\frac{n}{2(L^2 - 2LK)}\right]\right]$$

$$= \left(L\log(2p + \mathrm{e}^{m-1} - p) + \left[\frac{n}{2(L - 2K)}\right]\right)^m.$$

Now take

$$L = 2K + \sqrt{\frac{n}{2\log(2p + \mathrm{e}^{m-1} - p)}}\ .$$

□

**Lemma 8.1.** *(Jensen's inequality for partly concave functions) Let $X$ be a real-valued random variable, and let $g$ be an increasing function on $[0, \infty)$, which is concave on $[c, \infty)$ for some $c \geq 0$. Then*

$$\mathbf{E}g(|X|) \leq g\left[\mathbf{E}|X| + c\mathbf{P}(|X| < c)\right]\ .$$

**Proof.** We have

$$\mathbf{E}g(|X|) = \mathbf{E}g(|X|)\mathrm{l}\{|X| \geq c\} + \mathbf{E}g(|X|)\mathrm{l}\{|X| < c\}$$

$$\leq \mathbf{E}g(|X|)\mathrm{l}\{|X| \geq c\} + g(c)\mathbf{P}(|X| < c)$$

$$= \mathbf{E}\left[g(|X|)\Big||X| \geq c\right]\mathbf{P}(|X| \geq c) + g(c)\mathbf{P}(|X| < c) .$$

We now apply Jensen's inequality to the term on the left, and then use the concavity on $[c, \infty)$ to incorporate the term on the right:

$$\mathbf{E}g(|X|) \leq g\left[\mathbf{E}\left(|X|\Big||X| \geq c\right)\right]\mathbf{P}(|X| \geq c) + g(c)\mathbf{P}(|X| < c)$$

$$\leq g\left[\mathbf{E}|X| + c\mathbf{P}(|X| < c)\right] .$$

$\square$

### *8.2. Proofs for Section 3*

**Proof of Lemma 3.1.** This follows from

$$P(\gamma_f - \gamma_{f_0}) = P(l_f - l_{f_0}) .$$

$\square$

**Proof of Lemma 3.2.** We have

$$P|(f - f_0)(1 - 2\eta)| \geq vP|f - f_0|\mathrm{l}\{|1 - 2\eta| \geq v\}$$

$$\geq v\left(P|f - f_0| - P\mathrm{l}\{|1 - 2\eta| < v\}\right) := uv - H_1(v) ,$$

with $u = P|f - f_0|$. Since this is true for all $v$, we may maximize over $v$ to obtain

$$P|(f - f_0)(1 - 2\eta)| \geq G_1\left(P|f - f_0|\right) \geq G_1\left(P(f - f_0)^2\right) ,$$

as

$$P|f - f_0| \geq P(f - f_0)^2 .$$

Moreover,

$$|\gamma_f(y) - \gamma_{f_0}(y)| = |(f - f_0)(1 - 2y)| \leq |f - f_0| ,$$

so that

$$\sigma^2(\gamma_f - \gamma_{f_0}) \leq P(\gamma_f - \gamma_{f_0})^2 \leq P(f - f_0)^2 .$$

$\square$

**Proof of Lemma 3.3.** Clearly

$$P(\gamma_f - \gamma_{f_0}) = -\int_{f_0>0} \log\sqrt{\frac{f}{f_0}} f_0 d\mu$$

$$\geq -\int_{f_0>0} (\sqrt{\frac{f}{f_0}} - 1) f_0 d\mu$$

$$= 1 - \int \sqrt{f f_0} d\mu = h^2(f, f_0) \ .$$

Moreover,

$$\sigma^2(\gamma_{\bar{f}} - \gamma_{f_*}) \leq P(\gamma_{\bar{f}} - \gamma_{f_*})^2 \ .$$

Lemma 7.2 in (13) says that

$$2\{\exp|\gamma_{\bar{f}} - \gamma_{f_*}| - |\gamma_{\bar{f}} - \gamma_{f_*}| - 1\} \leq 8(\sqrt{\frac{\bar{f}}{f_*}} - 1)^2. \tag{13}$$

We moreover have

$$|\gamma_{\bar{f}} - \gamma_{f_*}|^2 \leq 2\{\exp|\gamma_{\bar{f}} - \gamma_{f_*}| - |\gamma_{\bar{f}} - \gamma_{f_*}| - 1\}$$

Thus

$$\sigma^2(\gamma_{\bar{f}} - \gamma_{f_*}) \leq 8 \int (\sqrt{\bar{f}} - \sqrt{f_*})^2 \frac{f_0}{f_*} d\mu \leq C^2 h^2(\bar{f}, f_*) \ .$$

$\square$

### *8.3. Proofs for Section 6*

**Proof of Lemma 6.1.** Define

$$\mathbf{Z} := \frac{|(P_n - P)(\hat{\gamma} - \gamma_*)|}{G^{-1}(\hat{\mathcal{E}}) + G^{-1}(\mathcal{E}_* \vee \varepsilon)} \ .$$

Then

$$\hat{\mathcal{E}} \leq \mathbf{Z}G^{-1}(\hat{\mathcal{E}}) + \mathbf{Z}G^{-1}(\mathcal{E}_* \vee \varepsilon) + \mathcal{E}_*$$

$$\leq \delta\hat{\mathcal{E}} + 2\delta H\left(\frac{\mathbf{Z}}{\delta}\right) + (1+\delta)\mathcal{E}_* \ .$$

It follows that

$$(1-\delta)\mathbf{E}\hat{\mathcal{E}} \leq 2\delta\mathbf{E}H\left(\frac{\mathbf{Z}}{\delta}\right) + (1+\delta)\mathcal{E}_*$$

$$\leq 2\delta H\left(\mathbf{E}\left(\frac{\mathbf{Z}}{\delta}\right)^r\right)^{1/r} + (1+\delta)\mathcal{E}_*$$

$$\leq 2\delta H\left(\sqrt{\frac{\Delta}{\delta^2}} + \frac{K\Delta}{\delta G^{-1}(\mathcal{E}_* \vee \varepsilon)}\right) + (1+\delta)\mathcal{E}_* \ .$$

$\square$

### *8.4. Proofs for Section 7*

#### *8.4.1. Preparatory lemmas*

We begin with two simple results (without proofs) for ease of reference.

**Lemma 8.2.** *If the loss envelope $\Gamma$ has power tails (2), then for all $m < 2s$ and $K > 0$,*

$$P\Gamma^{m/2}\mathbb{1}\{\Gamma > K\} \leq \frac{m}{2s - m}M^sK^{-(2s-m)/2} .$$

**Lemma 8.3.** *For positive constants $a, b, \alpha$ and $\beta$, the function*

$$g(x) := ax^{\alpha} + bx^{-\beta}, \; x > 0$$

*is minimized at*

$$x_0 := \left(\frac{b\beta}{a\alpha}\right)^{\frac{1}{\alpha+\beta}} ,$$

*and there attains a minimum of*

$$g(x_0) = \tilde{\mathcal{C}}(\alpha, \beta) \times a^{\frac{\beta}{\alpha+\beta}}b^{\frac{\alpha}{\alpha+\beta}} ,$$

*where*

$$\tilde{\mathcal{C}}(\alpha, \beta) := \left(\frac{\beta}{\alpha}\right)^{\frac{\alpha}{\alpha+\beta}} + \left(\frac{\alpha}{\beta}\right)^{\frac{\beta}{\alpha+\beta}} .$$

Next we need an auxiliary lemma:

**Lemma 8.4.** *For all $0 \leq z \leq 1$, we have that*

$$(1-z)^{2\kappa} \leq 1 - 2\kappa z^{2\kappa-1} + (2\kappa - 1)z^{2\kappa}$$

*and for all $z \geq 0$,*

$$(1+z)^{2\kappa} \geq 1 + 2\kappa z^{2\kappa-1} + z^{2\kappa} .$$

**Proof.** The second part is clear, as it involves the omission only of positive summands from the LHS to the RHS. For the first part, we write

$$f(z) := 1 - 2\kappa z^{2\kappa-1} + (2\kappa - 1)z^{2\kappa} - (1-z)^{2\kappa}$$

and note that

$$
\begin{aligned}
f(z) &= 1 - z^{2\kappa} - (1-z) \cdot 2\kappa z^{2\kappa-1} - (1-z)^{2\kappa} \\
&= (1-z) \cdot \left(\sum_{i=0}^{2\kappa-1} z^j - 2\kappa z^{2\kappa-1} - (1-z)^{2\kappa-1}\right) \\
&= (1-z)^2 \cdot \left(\sum_{j=0}^{2\kappa-2} (j+1) z^j - (1-z)^{2\kappa-2}\right) \\
&=: (1-z)^2 \cdot \tilde{f}(z) .
\end{aligned}
$$

Now as $\tilde{f}(0) = 0$ and for $0 \le z \le 1$,

$$\left(\tilde{f}\right)'(z) = \sum_{j=1}^{2\kappa-2} j(j+1)z^{j-1} + (2\kappa - 2) \cdot (1-z)^{2\kappa-3} \ge 0 \ ,$$

we know that $\tilde{f}(z)$, and thus $f(z)$, is non-negative on $[0,1]$ . $\qquad\square$

**Lemma 8.5.** *Let $a$, $b$ and $c$ be positive, let $\kappa \ge 1$, and assume that*

$$a \le b + c \cdot \left(a^{\frac{1}{2\kappa}} + b^{\frac{1}{2\kappa}}\right) \ .$$

*Then*

$$a^{\frac{1}{2\kappa}} \le \left(1 + (2\kappa-1)^{\frac{1}{2\kappa-1}}\right) \cdot \left(\frac{c}{2\kappa}\right)^{\frac{1}{2\kappa-1}} + b^{\frac{1}{2\kappa}} \ .$$

**Proof.** First note that if $a^{1/2\kappa} \le (c/2\kappa)^{1/(2\kappa-1)}$, then the desired inequality automatically holds. Thus we can restrict ourselves to the case where $a^{1/2\kappa} > (c/2\kappa)^{1/(2\kappa-1)}$. Applying the first part of Lemma 8.4 for $z = (c/2\kappa)^{1/(2\kappa-1)}/a^{1/2\kappa}$ – which now is less than 1 – gives us the inequality

$$\left(a^{\frac{1}{2\kappa}} - \left(\frac{c}{2\kappa}\right)^{\frac{1}{2\kappa-1}}\right)^{2\kappa} - \left(\frac{1}{2\kappa-1}\right)\left(\frac{c}{2\kappa}\right)^{\frac{2\kappa}{2\kappa-1}} \le a - c \cdot a^{\frac{1}{2\kappa}} \le b + c \cdot b^{\frac{1}{2\kappa}} \ ,$$

and thus

$$
\begin{aligned}
\left(a^{\frac{1}{2\kappa}} - \left(\frac{c}{2\kappa}\right)^{\frac{1}{2\kappa-1}}\right)^{2\kappa} &\le& b + c \cdot b^{\frac{1}{2\kappa}} + (2\kappa-1)\left(\frac{c}{2\kappa}\right)^{\frac{2\kappa}{2\kappa-1}} \\
&\le& b + (2\kappa-1)\,c \cdot b^{\frac{1}{2\kappa}} + \left(\frac{2\kappa-1}{2\kappa} \cdot c\right)^{\frac{2\kappa}{2\kappa-1}} \ ,
\end{aligned}
$$

where in the second step we used that $\kappa \ge 1$. Now part 2 of Lemma 8.4, applied to $z = \left(\frac{2\kappa-1}{2\kappa} \cdot c\right)^{1/2\kappa-1}/b^{1/2\kappa}$, yields

$$
\begin{aligned}
\left(b^{\frac{1}{2\kappa}} + \left(\frac{2\kappa-1}{2\kappa} \cdot c\right)^{\frac{1}{2\kappa-1}}\right)^{2\kappa} &\ge& b + (2\kappa-1) \cdot cb^{\frac{1}{2\kappa}} + \left(\frac{2\kappa-1}{2\kappa} \cdot c\right)^{\frac{2\kappa}{2\kappa-1}} \\
&\ge& \left(a^{\frac{1}{2\kappa}} - \left(\frac{c}{2\kappa}\right)^{\frac{1}{2\kappa-1}}\right)^{2\kappa} \ ,
\end{aligned}
$$

from which the stated inequality follows. $\qquad\square$

*8.4.2. Main proof*

**Proof of Theorem 7.1.** (i) In the power tail case, we define

$$\mathcal{E}_*^\tau := \mathcal{E}_* \vee \tau \ ,$$

where $\tau$ is a strictly positive number, and

$$\mathbf{Z} := \frac{\left| P_n \left( (\hat{\gamma}^c - \gamma_*^c) \, 1 \, \{\Gamma \le K\} \right)^c \right|}{C \left( \hat{\mathcal{E}}^{\frac{1}{2\kappa}} + (\mathcal{E}_*^\tau)^{\frac{1}{2\kappa}} \right)} \ .$$

Then we have

$$
\begin{aligned}
\hat{\mathcal{E}} &\le \ |(P_n - P)(\hat{\gamma} - \gamma_*)| + \mathcal{E}_* \\
&= \ |P_n \left( \hat{\gamma}^c - \gamma_*^c \right)| + \mathcal{E}_* \\
&\le \ |P_n \left( (\hat{\gamma}^c - \gamma_*^c) \, 1 \, \{\Gamma \le K\} \right)^c| + |P \left( (\hat{\gamma}^c - \gamma_*^c) 1 \, \{\Gamma \le K\} \right)| \\
&\quad + |P_n \left( (\hat{\gamma}^c - \gamma_*^c) \, 1 \, \{\Gamma > K\} \right)| + \mathcal{E}_* \\
&\le \ C \mathbf{Z} \left( \hat{\mathcal{E}}^{\frac{1}{2\kappa}} + (\mathcal{E}_*^\tau)^{\frac{1}{2\kappa}} \right) + \mathcal{E}_* + (P_n + P) \left( \Gamma 1 \, \{\Gamma > K\} \right) \\
&\le \ C \mathbf{Z} \left( \hat{\mathcal{E}}^{\frac{1}{2\kappa}} + (\mathcal{E}_*^\tau + (P_n + P) \left( \Gamma 1 \, \{\Gamma > K\} \right))^{\frac{1}{2\kappa}} \right) \\
&\quad + \mathcal{E}_*^\tau + (P_n + P) \left( \Gamma 1 \, \{\Gamma > K\} \right) \ .
\end{aligned}
$$

Using Lemma 8.5, we obtain the inequality

$$
\begin{aligned}
\hat{\mathcal{E}}^{\frac{1}{2\kappa}} &\le \ \left( 1 + (2\kappa - 1)^{\frac{1}{2\kappa - 1}} \right) \left( \frac{C \mathbf{Z}}{2\kappa} \right)^{\frac{1}{2\kappa - 1}} \\
&\quad + (\mathcal{E}_*^\tau + (P_n + P) \left( \Gamma 1 \, \{\Gamma > K\} \right))^{\frac{1}{2\kappa}} \\
&\le \ \left( 1 + (2\kappa - 1)^{\frac{1}{2\kappa - 1}} \right) \left( \frac{C \mathbf{Z}}{2\kappa} \right)^{\frac{1}{2\kappa - 1}} \\
&\quad + (\mathcal{E}_*^\tau)^{\frac{1}{2\kappa}} + ((P_n + P) \left( \Gamma 1 \, \{\Gamma > K\} \right))^{\frac{1}{2\kappa}} \ ,
\end{aligned}
$$

where for the second step we used the elementary observation $a^{2\kappa} + b^{2\kappa} \le (a + b)^{2\kappa}$ for $a, b \ge 0$, $\kappa > 1/2$. Now we will first compute the moments of $\mathbf{Z}$ by an application of Bernstein's inequality. We know that

$$ P \left| (\gamma_j^c - \gamma_*^c) \, 1 \, \{\Gamma \le K\} \right|^m \ \le \ K^{m-2} P \left[ \left( (\gamma_j^c - \gamma_*^c) \, 1 \, \{\Gamma \le K\} \right)^2 \right] $$

and

$$
\begin{aligned}
P \left[ \left( (\gamma_j^c - \gamma_*^c) \, 1 \, \{\Gamma \le K\} \right)^2 \right] &= \ P \left[ (\gamma_j^c - \gamma_*^c)^2 \, 1 \, \{\Gamma \le K\} \right] \\
&\le \ P \left[ (\gamma_j^c - \gamma_*^c)^2 \right] \\
&= \ \sigma^2 \left( \gamma_j - \gamma_* \right) \\
&= \ \sigma^2 \left( (\gamma_j - \gamma_0) - (\gamma_* - \gamma_0) \right) \\
&\le \ \left( \sigma \left( \gamma_j - \gamma_0 \right) + \sigma \left( \gamma_* - \gamma_0 \right) \right)^2 \ ,
\end{aligned}
$$

which by the margin condition

$$
\begin{aligned}
&\le \ \left( C \cdot \left( P \left( \gamma_j - \gamma_0 \right) \right)^{1/2\kappa} + C \cdot \left( P \left( \gamma_* - \gamma_0 \right) \right)^{1/2\kappa} \right)^2 \\
&= \ C^2 \cdot \left( \mathcal{E}_j^{1/2\kappa} + \mathcal{E}_*^{1/2\kappa} \right)^2 \ .
\end{aligned}
$$

Thus for all $j$,

$$P \left| \frac{(\gamma_j^c - \gamma_*^c) \, 1\{\Gamma \leq K\}}{C \left(\mathcal{E}_j^{1/2\kappa} + (\mathcal{E}_*^\tau)^{1/2\kappa}\right)} \right|^m \quad \leq \quad \left( \frac{K}{C \left(\mathcal{E}_j^{1/2\kappa} + (\mathcal{E}_*^\tau)^{1/2\kappa}\right)} \right)^{m-2}$$

$$\leq \quad \left( \frac{K}{C \, (\mathcal{E}_*^\tau)^{1/2\kappa}} \right)^{m-2}$$

$$\Rightarrow P \left| \frac{\left((\gamma_j^c - \gamma_*^c) \, 1\{\Gamma \leq K\}\right)^c}{C \left(\mathcal{E}_j^{1/2\kappa} + (\mathcal{E}_*^\tau)^{1/2\kappa}\right)} \right|^m \quad \leq \quad 4 \cdot \left( \frac{2K}{C \, (\mathcal{E}_*^\tau)^{1/2\kappa}} \right)^{m-2} \;,$$

and we can apply Corollary 2.1 to obtain

$$\|\mathbf{Z}\|_m = \left\| \frac{P_n \left[(\hat{\gamma}^c - \gamma_*^c) \, 1\{\Gamma \leq K\}\right]}{C \left(\hat{\mathcal{E}}^{1/2\kappa} + (\mathcal{E}_*^\tau)^{1/2\kappa}\right)} \right\|_m \quad \leq \quad 2 \left( \sqrt{\Delta} + \frac{K\Delta}{C \, (\mathcal{E}_*^\tau)^{1/2\kappa}} \right) \;.$$

Now to compute the moments of

$$(P_n + P) \left(\Gamma 1\{\Gamma > K\}\right)^{\frac{1}{2\kappa}} \;,$$

we proceed as follows for $m \geq 2\kappa$ (using that $\kappa \geq 1/2$):

$$\left\| \left((P_n + P)\left(\Gamma 1\{\Gamma > K\}\right)\right)^{1/2\kappa} \right\|_m$$

$$= \quad \left( \mathbf{E}\left[ \left((P_n + P)\left(\Gamma 1\{\Gamma > K\}\right)\right)^{m/2\kappa} \right] \right)^{1/m}$$

$$\leq \quad \left( 2^{m/2\kappa - 1} \mathbf{E}\left[ (P_n + P)\left(\Gamma^{m/2\kappa} 1\{\Gamma > K\}\right) \right] \right)^{1/m}$$

$$= \quad 2^{1/2\kappa} \left( P\left(\Gamma^{m/2\kappa} 1\{\Gamma > K\}\right) \right)^{1/m} \;.$$

By Lemma 8.2, for $m < 2s\kappa$, this has an upper bound in

$$2^{1/2\kappa} \left( \frac{m}{2s\kappa - m} \right)^{1/m} M^{s/m} K^{1/2\kappa - s/m} \;.$$

Thus we find that for $m \in [2\kappa, \min\{1 + \log(p), 2s\kappa\})$ (and remembering that $\tau^2 = \mathcal{E}$),

$$\left\| \left(\hat{\mathcal{E}}\right)^{\frac{1}{2\kappa}} \right\|_m \leq (\mathcal{E}_*^\tau)^{\frac{1}{2\kappa}} + A(\kappa) \cdot C^{\frac{1}{2\kappa - 1}} \cdot \left( \sqrt{\Delta} + \frac{K\Delta}{C(\mathcal{E}_*^\tau)^{\frac{1}{2\kappa}}} \right)^{\frac{1}{2\kappa - 1}}$$

$$+ B(\kappa, s, m) \cdot M^{s/m} K^{1/2\kappa - s/m} \;,$$

where

$$A(\kappa) := \frac{1 + (2\kappa - 1)^{\frac{1}{2\kappa - 1}}}{\kappa^{\frac{1}{2\kappa - 1}}} \;,$$

$$B(\kappa, s, m) := 2^{1/2\kappa} \left( \frac{m}{2s\kappa - m} \right)^{\frac{1}{m}}.$$

If we now apply the straightforward bound

$$\left( \sqrt{\Delta} + \frac{K\Delta}{C(\mathcal{E}_*^\tau)^{\frac{1}{2\kappa}}} \right)^{\frac{1}{2\kappa-1}} \leq \left( \sqrt{\Delta} \right)^{\frac{1}{2\kappa-1}} + \left( \frac{K\Delta}{C(\mathcal{E}_*^\tau)^{\frac{1}{2\kappa}}} \right)^{\frac{1}{2\kappa-1}}$$

and minimize the upper bound over $K \geq 0$ (using Lemma 8.3), we obtain the desired oracle inequality for the power tail case.

(ii) If we assume the exponential moment condition instead of power tails, we can take

$$\mathbf{Z} := \frac{|P_n\left( (\hat{\gamma}^c - \gamma_*^c) \right)|}{C \left( \hat{\mathcal{E}}^{\frac{1}{2\kappa}} + (\mathcal{E}_*^\tau)^{\frac{1}{2\kappa}} \right)}$$

and we obtain the same bound for $\|\mathbf{Z}\|_m$ as before, but no term stemming from $\Gamma 1 \{\Gamma > K\}$. This yields the desired risk moment inequality. The corresponding risk tail bound also comes straight from applying Bernstein's inequality (3) to $\mathbf{Z}$.

$$\square$$